\documentstyle[12pt,graphicx]{article}
\def\thebibliography#1
{\begin{center}{\normalsize\bf References}\end{center}
\list{[\arabic{enumi}]}%
{\settowidth\labelwidth{[#1]}\leftmargin \labelwidth
\advance\leftmargin\labelsep
\usecounter{enumi}}
\sloppy\clubpenalty4000\widowpenalty4000}

\begin{document}
\baselineskip=16pt

\title{\vspace*{0cm}\large BRANCHED COVERS OF TANGLES IN THREE-BALLS}

\author{{\normalsize MAKIKO ISHIWATA}\\
{\small Department of Mathematics, Tokyo Woman's Christian University}\\[-1mm]
{\small Zempukuji 2-6-1, Suginamiku, Tokyo 167-8585, Japan}\\[-1mm]
{\small e-mail: mako@twcu.ac.jp}
\\[3mm]
{\normalsize J\'{O}ZEF H. PRZYTYCKI}\\
{\small Department of Mathematics, The George Washington University}\\[-1mm]
{\small Washington, DC 20052, USA}\\[-1mm]
{\small e-mail: przytyck@research.circ.gwu.edu}
\\[3mm]
{\normalsize AKIRA YASUHARA }\\
{\small Department of Mathematics, Tokyo Gakugei University}\\[-1mm]
{\small Nukuikita 4-1-1, Koganei, Tokyo 184-8501, Japan}\\[-1mm]
{\small e-mail: yasuhara@u-gakugei.ac.jp}
}

\date{\empty}

\maketitle

\begin{abstract}
We give an algorithm for a surgery description of a $p$-fold
 cyclic branched cover of $B^3$ branched along a tangle. We generalize
constructions of Montesinos and Akbulut-Kirby.
\end{abstract}



\vspace*{10mm}
Tangles were first studied by Conway \cite{Con}. They were particularly
useful for analyzing prime and hyperblic knots.
A branched cover of a tangle is an indispensable tool for understanding 
tangles. Hence it is important to give practical presentations of
branched covers of tangles.
A $p$-fold branched cover of an $n$-tangle is a three-manifold
with the boundary a connected surface of genus $(n-1)(p-1)$. Such a manifold 
can be obtained from the genus $(n-1)(p-1)$ handlebody by a surgery.
We provide an algorithm for a surgery description of a $p$-fold cyclic
branched cover of $B^3$ branched along a tangle.
 The construction generalizes that of Montesinos \cite{Mon} and 
 Akbulut and Kirby \cite{A-K}. 
It is strikingly simple in the case of the two-fold branched cover.
We also discuss the related Heegaard decomposition of a $p$-fold
branched cover of an $n$-tangle.

\bigskip
\noindent
{\bf{1. Sugery descriptions}}

\medskip
A ${\it{tangle}}$ is a one-manifold properly embedded in a
three-ball. An $n$-tangle is a tangle with $2n$ boundary points.
Let $T$ be an $n$-tangle and $T_0$ a trivial 
$n$-tangle diagram\footnote{Tangles are considered up to ambient 
isotopy but in practice
we will often use a tangle for a tangle diagram or an actual embedding 
of a one-manifold.} (Fig.~1).
Let $D_1 \cup \cdots \cup D_n$ be a 
disjoint union of disks bounded by
$T_0$ and $b_1, \cdots, b_m$ be mutually disjoint disks
in $B^3$ such that $b_i \cap {\bigcup}_j D_j=\partial b_i \cap T_0$ are
two 
disjoint arcs in $\partial b_i (i=1,\cdots,m)$, see Fig.~2. 
We denote by $\Omega (T_0; \{D_1,\cdots,D_n\},\{b_1,
 \cdots,b_m\}) $ the tangle
 $T_0 \cup {\bigcup}_i \partial b_i - \mbox{int}(T_0 \cap 
 {\bigcup}_i\partial b_i )$ 
and call it a {\em disk-band representation}\/\footnote{It
contains information on the tangle, its surface and the decomposition of 
the surface.} of a tangle. 
A disk-band representation is called {\em bicollared} 
if the surface ${\bigcup}_i D_i \cup {\bigcup}_j b_j$ is orientable.
We will see that any $n$-tangle has a bicollared 
disk-band representation (Proposition~5).

\begin{center}
\begin{tabular}{cc} 
\includegraphics[trim=0mm 0mm 0mm 0mm, width=.25\linewidth]
{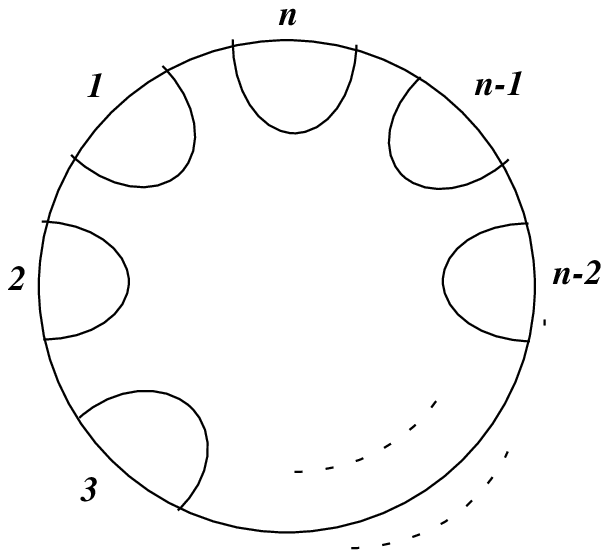}
&\hspace*{20mm}
\includegraphics[trim=0mm 0mm 0mm 0mm, width=.25\linewidth]
{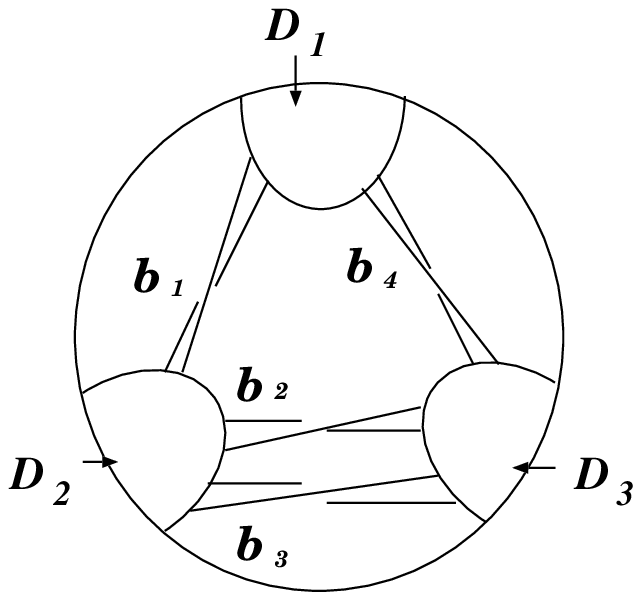}
\\
Fig. 1 &\hspace*{20mm} Fig. 2
\end{tabular}
\end{center}

\medskip
A {\em framed link} is a disjoint union of embedded annuli in a 3-manifold.
Framed links in $S^3$ can be identified with links whose each component 
is assigned an integer. Such links are also called framed links. 
Let $M$ be a three-manifold and $\cal{L}$ a framed link in $M$.
We denote by $\Sigma({\cal{L}},M)$ the manifold that is obtained from
$M$ by the surgery along $\cal{L}$ \cite{P-S}.

The case of two-fold branched covers is easy to visualize 
so we formulate it first.

\medskip
\noindent
{\bf{Theorem 1.}}{\em{
Let $\Omega (T_0; \{D_1,\cdots,D_n\},\{b_1,\cdots,b_m\})$ 
be a disk-band representation of an $n$-tangle $T$ in $B^3$.
Let $\varphi : H_0 \rightarrow B^3 $ be the two-fold branched cover of
 $B^3$ by a genus $n-1$ handlebody $H_0$ 
 branched along $T_0$. Then the two-fold branched cover
 of $B^3$ branched along $T$ has a surgery description
 $\Sigma({\varphi}^{-1}(\bigcup_i b_i),H_0)$ (compare Fig.~3).}}

\begin{center}
\begin{tabular}{c} 
\includegraphics[trim=0mm 0mm 0mm 0mm, width=.65\linewidth]
{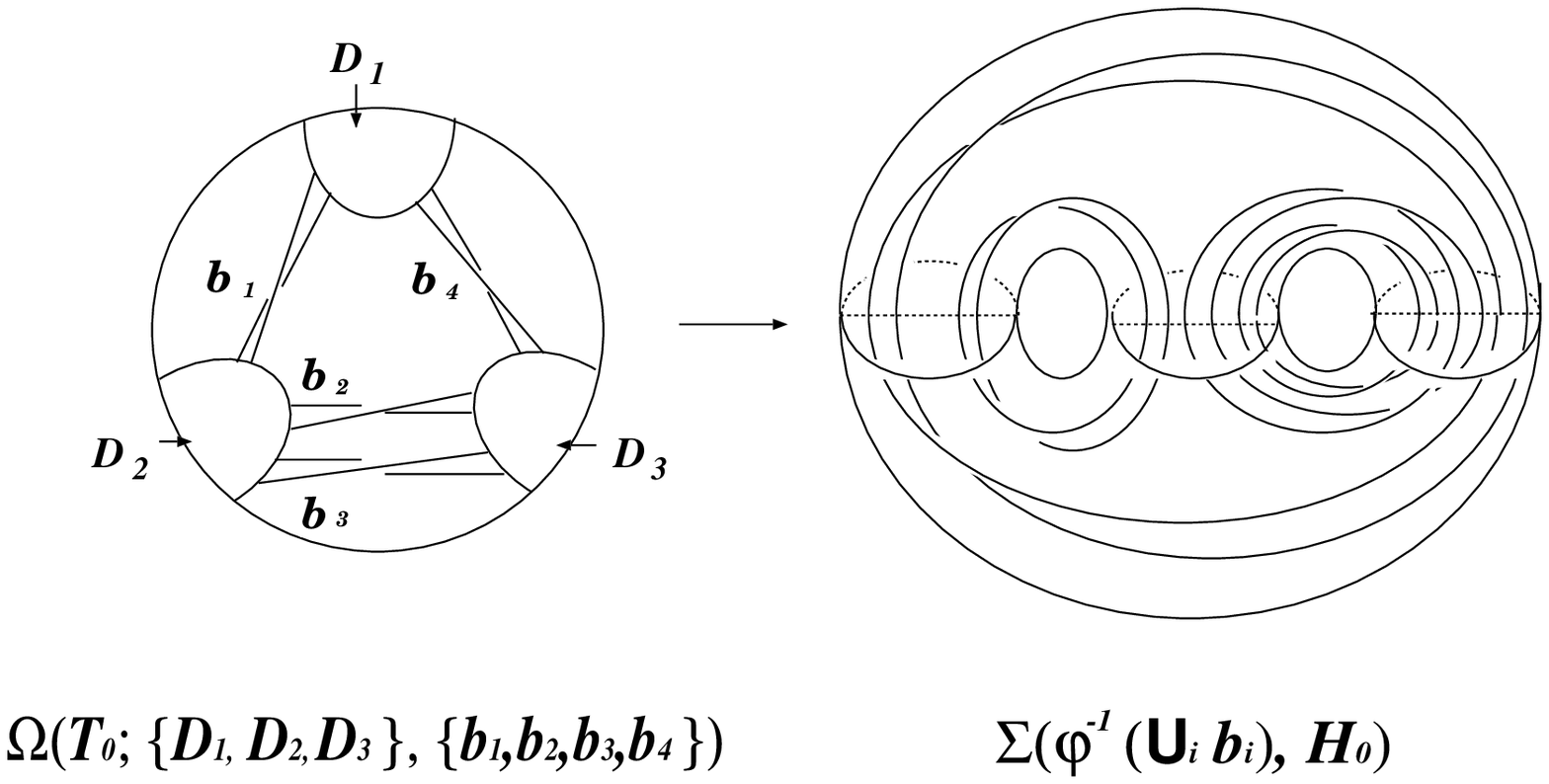}
\\
Fig. 3
\end{tabular}
\end{center}

\medskip
{\it{Proof.}} Let $X$ be $B^3-{\bigcup}_i {D_i}$ compactified with 
two copies, $D_i^{\pm}$,
of $D_i~~(i=1,2,\cdots,n)$ (Fig.~4). Let $X_1$ and
$X_2$ be two copies of $X$, and let $D_{i,k}^{\pm} \subset X_k$ 
denote copies of $D_{i}^{\pm} ~(i=1,2,\cdots,n,~ k = 1,2)$ (Fig.~5).
Then $H_0$ is obtained from $X_1 \cup X_2$ by
 identifying $D_{i,1}^{\varepsilon}$ with
 $D_{i,2}^{- \varepsilon}$~(${\varepsilon} \in \{-,+\}$). 
Let $b_{j,k}={\varphi}^{-1}(b_j) \cap X_k$ and let $Y$ be 
$H_0 - {\bigcup}_{j,k} b_{j,k}$ compactified with two copies
$b_{j,k}^{\pm}$ of $b_{j,k}$~in $X_k~(j=1,2,\cdots,m,~ k=1,2)$.
Here, $+$ or $-$ sides of $D_{i,k}$ and $b_{j,k}$ are not necessarily
compatible. We note, and it is the key observation of the construction,
that the two-fold branched cover $H$ of $B^3$ branched
along $T$ is obtained from $Y$ by identifying $b_{j,1}^{\varepsilon}$
with $b_{j,2}^{- \varepsilon}~(\varepsilon \in \{-,+\})$.
Note that each $b_{j,1}^+ \cup b_{j,2}^-\cup b_{j,1}^- \cup b_{j,2}^+$ 
is a torus. Let $c_j$ be the core of the annulus 
$b_{j,1}^+ \cup b_{j,2}^-$. The manifold obtained from $Y$ 
by identifying $b_{j,1}^{\varepsilon}$ with $b_{j,2}^{- \varepsilon}$ 
is homeomorphic to the one obtained from $Y$ by 
attaching tori $D_j^2\times S^1~(j=1,2,...,m)$ so that 
$\partial D_j^2=c_j$. Hence $H$ is homeomorphic to the manifold with the
surgery description 
$\Sigma({\varphi}^{-1}({\bigcup}_j b_j), H_0)$.
\hfill$\Box$

\begin{center}
\begin{tabular}{cc} 
\includegraphics[trim=0mm 0mm 0mm 0mm, width=.45\linewidth]
{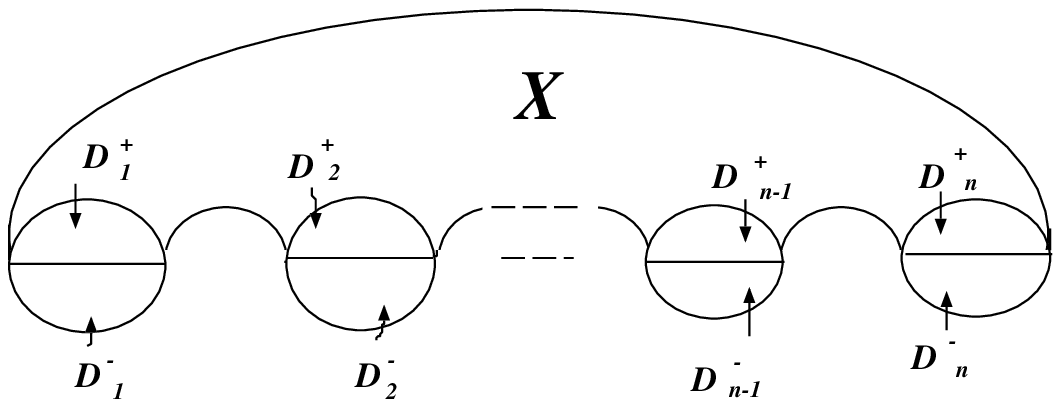}
&\hspace*{10mm}
\includegraphics[trim=0mm 0mm 0mm 0mm, width=.45\linewidth]
{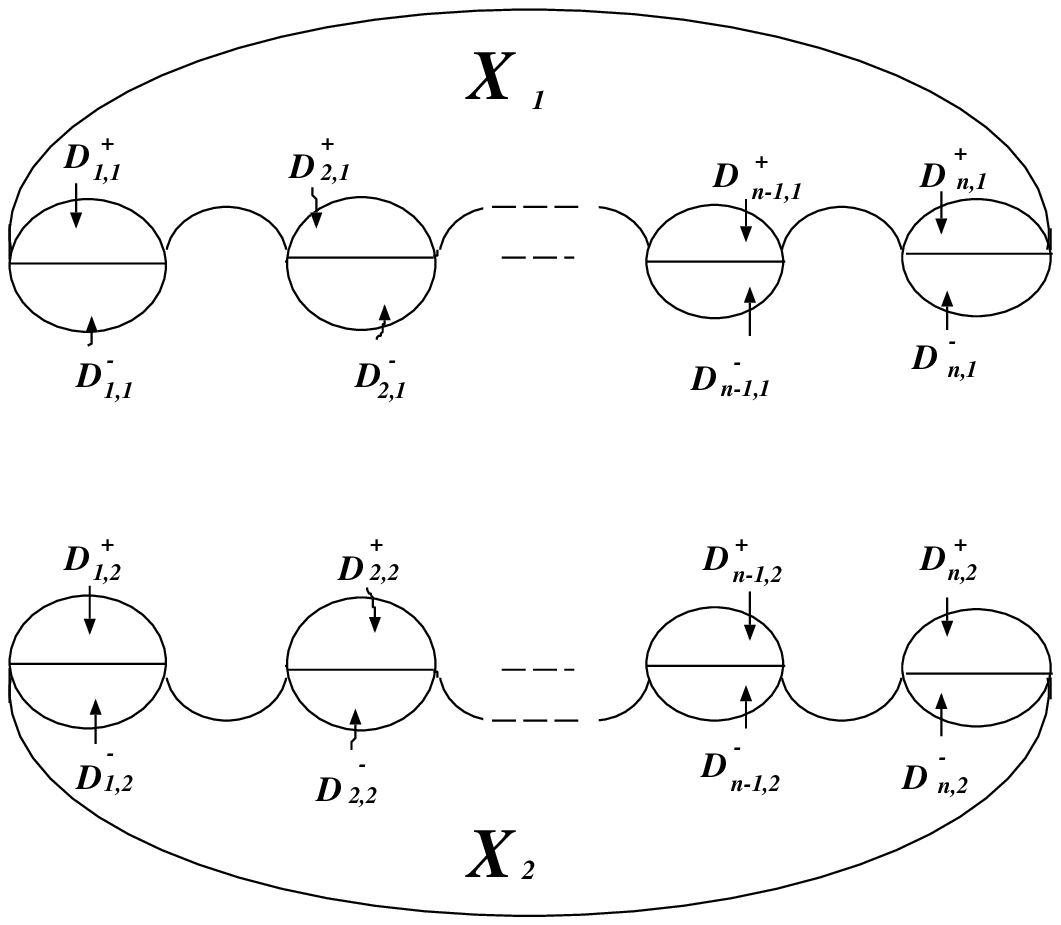}
\\
Fig. 4 &\hspace*{10mm} Fig. 5
\end{tabular}
\end{center}

\medskip
\noindent
{\bf{Example 2.}} (a) The two-fold branched cover $M^{(2)}(T_1)$ branched 
along a tangle $T_1$ in Fig.~6
is the Seifert manifold with the base a disk and two special fibers of
type $(2,1)$ and $(2,-1)$. Furthermore $M^{(2)}(T_1)$ is a twisted $I$-bundle
over the Klein bottle (for example see \cite{Jac}).
In particular, ${\pi}_1(M^{(2)}(T_1))=\langle a,b \mid aba^{-1}b=1 \rangle$. 

\noindent
(b) If we glue together two copies of $T_1$ as in Fig.~7, we get a
Borromean rings $L$. Thus our previous computation shows that the
two-fold branched cover $M^{(2)}(L)$ of $S^3$ branched along $L$ is a
``switched'' double of the twisted $I$-bundle over the
Klein bottle (see Fig.~8 for a surgery description). The fundamental
group ${\pi}_1(M^{(2)}(L))=
\langle x,a \mid x^2ax^2a^{-1}, a^2xa^2x^{-1}\rangle$
is a $3$-manifold group which is torsion 
free but not left orderable \cite{Rol}.

\noindent
(c) If we take the double of the tangle $T_1$, we obtain the link in
Fig.~9. The two-fold branched cover of $S^3$ branched along this link
is the double of twisted $I$-bundle over Klein bottle. 
A surgery description of this manifold is shown in Fig.~10. 
Thus this manifold is the Seifert manifold of type
$(2,1),(2,1),(2,-1),(2,-1)$. This manifold also has another
Seifert fibration, which is a circle bundle over 
the Klein bottle.

\medskip
\noindent
Example 2(a) was motivated by the fact that the tangle $T_1$ yields a 
virtual Lagrangian of index $2$ in the symplectic space of the Fox 
${\bf Z}$-colorings of the boundary of our tangle \cite{DJP}.

\begin{center}
\begin{tabular}{c} 
\includegraphics[trim=0mm 0mm 0mm 0mm, width=.5\linewidth]
{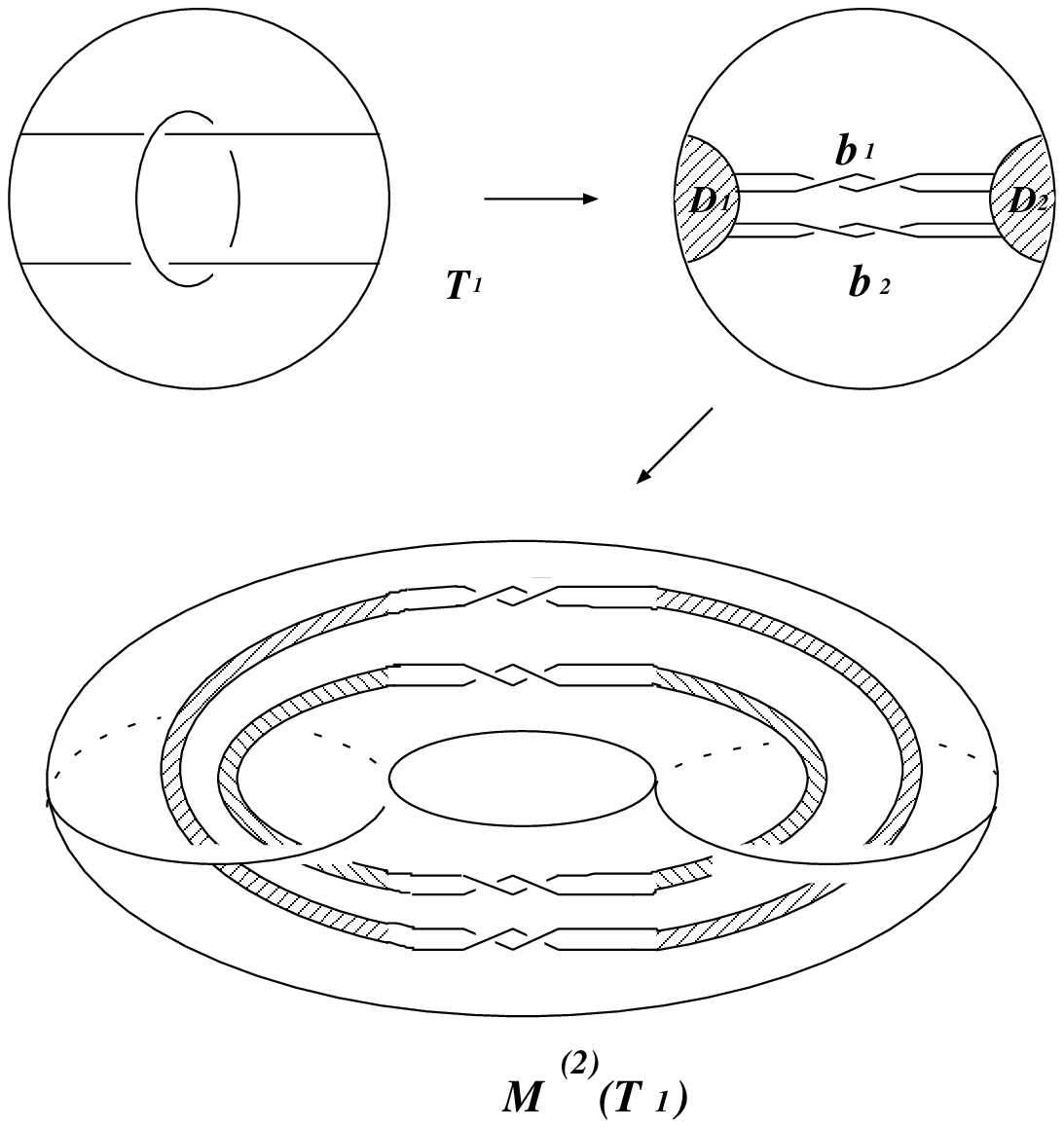}
\\
Fig. 6 
\end{tabular}

\begin{tabular}{cccc} 
\includegraphics[trim=0mm 0mm 0mm 0mm, width=.18\linewidth]
{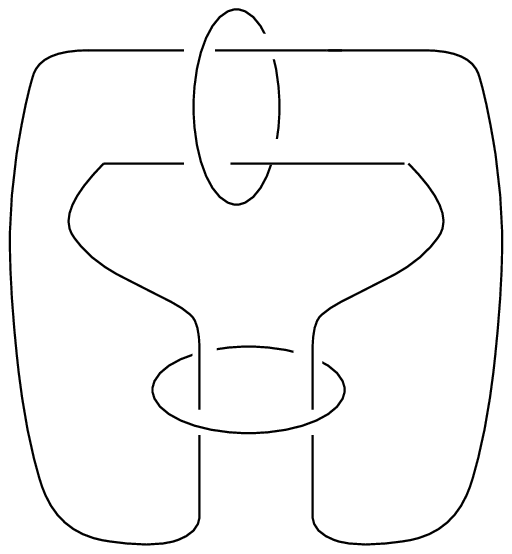}
&
\includegraphics[trim=0mm 0mm 0mm 0mm, width=.27\linewidth]
{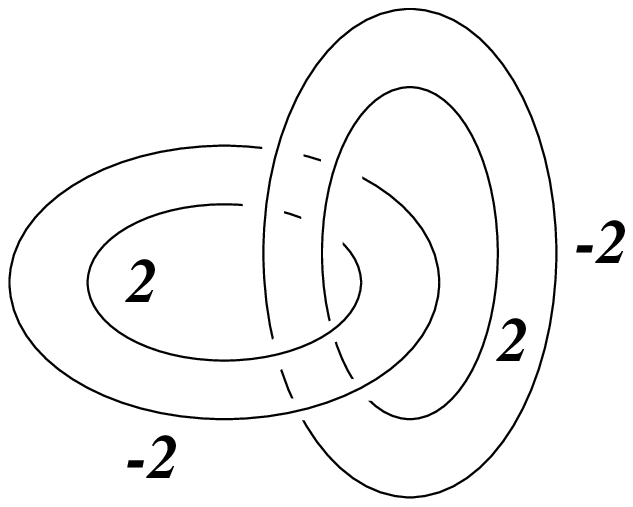}
&
\includegraphics[trim=0mm 0mm 0mm 0mm, width=.2\linewidth]
{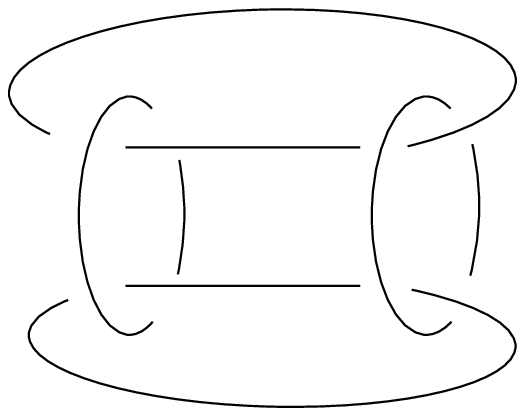}
&
\includegraphics[trim=0mm 0mm 0mm 0mm, width=.27\linewidth]
{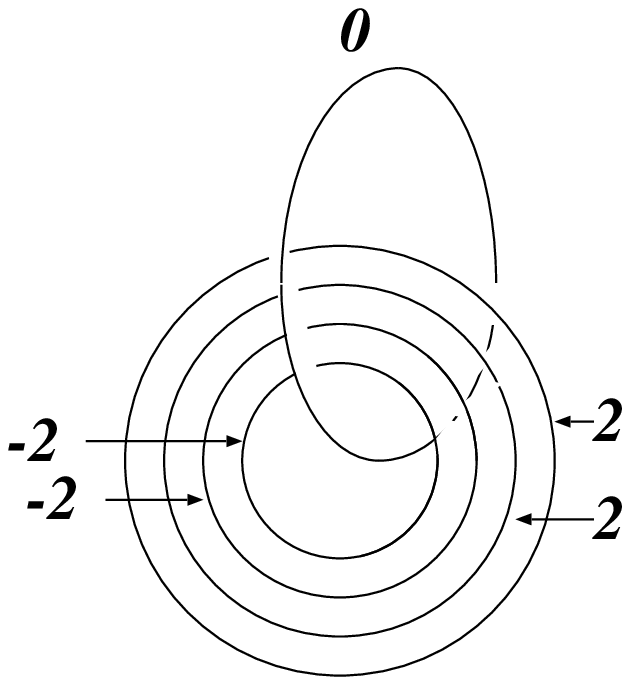}
\\
Fig. 7 & Fig. 8 & Fig. 9 & Fig. 10
\end{tabular}


\end{center}

\medskip
More generally we have:

\medskip
\noindent
{\bf{Example 3.}} Consider a tangle $T_2$ in Fig.~11, called a
{\it{pretzel tangle}} of type 
$(a_1, a_2, \cdots, a_n)$, where each $a_i$ is an integer 
indicating the number of half-twists $(i=1,2,\cdots,n)$.
The two-fold branched cover $M^{(2)}(T_2)$ branched along the tangle $T_2$ is a
 Seifert fibered manifold with the base a disk and $n$
 special fibers of type $(a_1,1)$,$(a_2,1), \cdots, (a_n,1)$ (Fig.~12).

\begin{center}
\begin{tabular}{cc} 
\includegraphics[trim=0mm 0mm 0mm 0mm, width=.3\linewidth]
{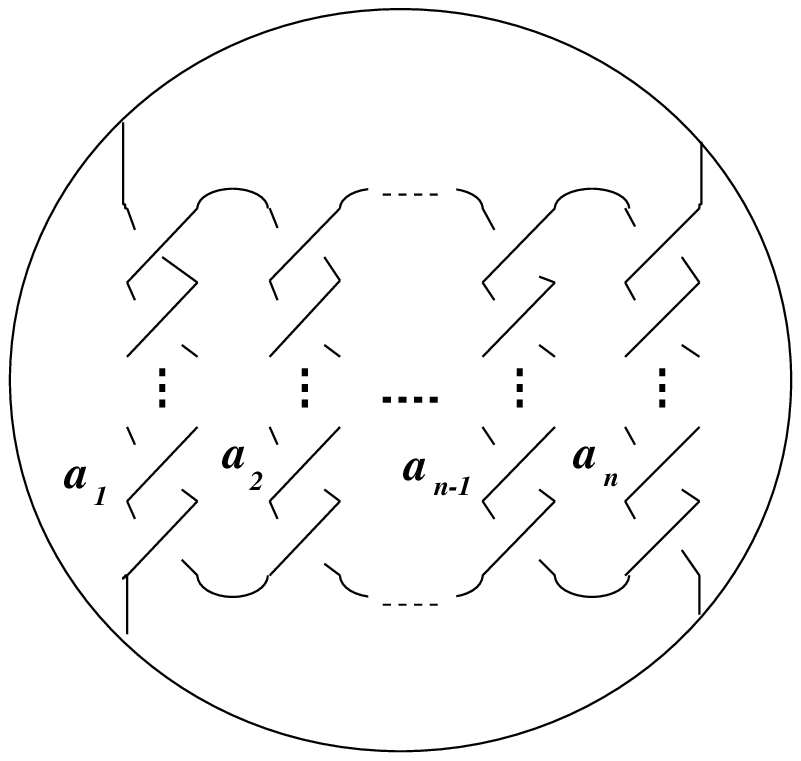}
&\hspace*{10mm}
\includegraphics[trim=0mm 0mm 0mm 0mm, width=.5\linewidth]
{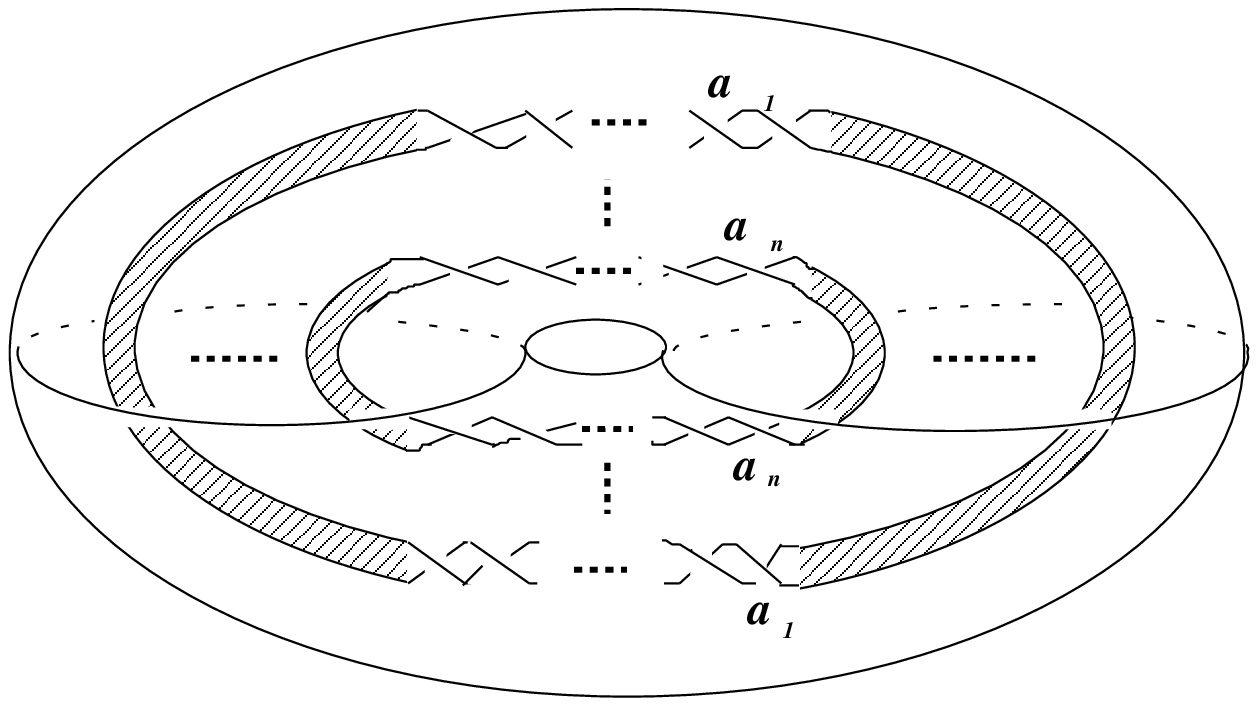}
\\
Fig. 11 &\hspace*{10mm} Fig. 12
\end{tabular}
\end{center}

\medskip
Theorem 1 can be generalized to the $p$-fold cyclic branched
cover assuming that a disk-band representation is bicollared.
We proceed as follows:

\noindent
Let $T=\Omega (T_0; \{D_1,\cdots,D_n\},\{b_1,\cdots,b_m\})$
 be a 
bicollared disk-band representation of  an $n$-tangle. 
Then ${\bigcup}_i D_i \cup {\bigcup}_j b_j$ has 
bicollar neighbourhood 
$({\bigcup}_i D_i \cup {\bigcup}_j b_j)\times [-1,1]$.
Let $X=\overline{B^3-((\bigcup_i {D_i}) \times [-1,1])}$ 
and ${D_i}^{\pm}=(D_i \times
[\pm 1,0]) \cap \partial X$. Let $X_k$ be a copy of $X$ and $D_{i,k}^{\pm} 
\subset \partial X_k$ a copy of $D_i^{\pm}  (k =1,2,\cdots,p)$. Then the
$p$-fold cyclic branched cover ${\varphi} : H_0 \rightarrow B^3$ 
branched along
$T_0$ is obtained from $X_1 \cup \cdots \cup X_p$ by identifying
$D_{i,k}^+$ with  $D_{i,k+1}^- (k=1,\cdots,p)$, where $k$ is considered
modulo $p$. Let $b_{j,k}^{\pm}={\varphi}^{-1}(b_j \times \{\pm 1\}) \cap X_k$. 
Note that each $b_{j,k}^+ \cup b_{j,k+1}^-$ is an annulus in $H_0$ 
for any $j$  and $k$. 
Then we obtain the $p$-fold cyclic branched cover of $B^3$
branched along $T$ in the similar way as in Theorem~1.

\medskip
\noindent
{\bf{Theorem 4.}}
{\em{Let $\Omega (T_0; \{D_1,\cdots,D_n\},\{b_1,\cdots,b_m\})$ be a 
bicollared disk-band representation of  an $n$-tangle $T$ in $B^3$.
Then 
$\Sigma({\bigcup}_{j=1}^m( {\bigcup}_{k=1}^{p-1} (b_{j,k}^+ \cup
 b_{j,k+1}^-)),H_0)$  is the $p$-fold cyclic branched cover
 of $B^3$ branched along $T$.}}  

\medskip
\noindent
Note that we do not use the annuli $b_{j,p+1}^+ \cup b_{j,1}^-$ 
$(j=1,2,\cdots,m)$ in the theorem above. In fact the cores of these
annuli bound mutually disjoint 2-dsks in 
$\Sigma({\bigcup}_{j=1}^m( {\bigcup}_{k=1}^{p-1} (b_{j,k}^+ \cup
 b_{j,k+1}^-)),H_0)$. 

\medskip
{\it{Proof.}} Let $Y=\overline{H_0 - {\bigcup}_j {\varphi}^{-1}
(b_j \times [-1,1])}$,
$V_{j,k}^{\pm}={\varphi}^{-1}(b_j \times [\pm 1,0]) \cap X_k$ and
${\beta}_{j,k}^{\pm} = V_{j,k}^{\pm} \cap Y (= \partial V_{j,k}^{\pm} \cap
\partial Y )$. Note that ${\bigcup}_j {\varphi}^{-1}(b_j \times
[-1,1])$ is a genus $p-1$ handlebody. Then the $p$-fold cyclic branched 
cover of $B^3$ branched along $T$ is
homeomorphic to a manifold $H$  that is obtained from $Y$ by identifying 
${\beta}_{j,k}^+$ with ${\beta}_{j,k+1}^- (k=1,\cdots,p)$, where
$k$ is taken modulo $p$. Moreover 
$H$ is homeomorphic to a manifold obtained from $\overline{H_0
- {\bigcup}_j {\varphi}^{-1}(b_j \times [-1,1])} \cup {\bigcup}_j 
(V_{j,1}^- \cup
V_{j,p}^+ \cup {\varphi}^{-1}(b_j \times \{0\})) $ by identifying
${\beta}_{j,k}^+$ and $b_{j,k}$ with ${\beta}_{j,k+1}^-$ and $b_{j,k+1}
(j=1,\cdots,m,~k=1,\cdots,p-1)$ respectively, where
$b_{j,k}={\varphi}^{-1}(b_j \times \{0\}) \cap X_k$. Note that
${\beta}_{j,k}^+ \cup b_{j,k} \cup {\beta}_{j,k+1}^- \cup b_{j,k+1}$ is a
torus. By the argument similar to that in the proof of Theorem 1, we
have the required result.  
\hfill$\Box$

\medskip
\noindent
{\bf{Example 5.}} The three-fold cyclic branched cover $M^{(3)}(4_1)$ 
of $S^3$ branched along the figure eight knot has a surgery description 
shown in Fig.~13(a). The framed link in Fig. 13(a) can be 
deformed into the link in Fig.~13(b) by an ambient isotopy and 
a second Kirby move. 
The link in Fig~13(b) is ambient isotopic to the link in Fig~8. 
Hence $M^{(3)}(4_1)$ is homeomorphic to the two-fold branched cover of 
$S^3$ branched along the Borromean rings (cf. Example 2(b)).

\begin{center}
\begin{tabular}{c} 
\includegraphics[trim=0mm 0mm 0mm 0mm, width=.95\linewidth]
{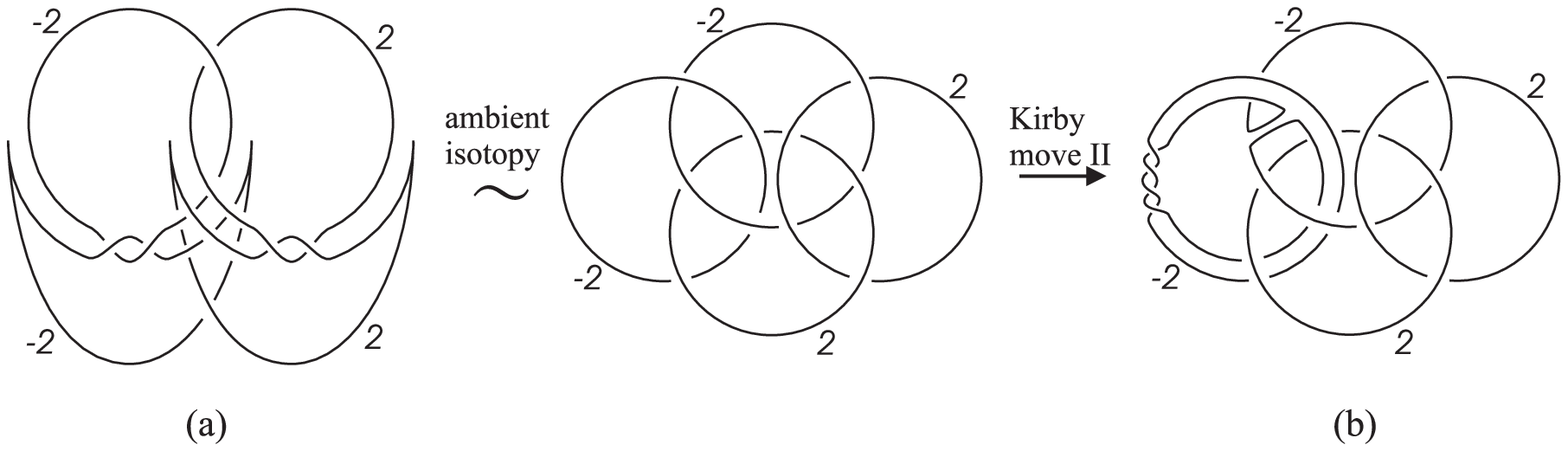}
\\
Fig. 13
\end{tabular}
\end{center}

\medskip
\noindent
{\bf{Proposition 6.}}  
{\em{Any $n$-tangle $(B^3,T)$ has a (bicollared) disk-band representation.}} 

\medskip 
{\it Proof.} 
We attach a trivial $n$-tangle $T_0$ to
the $n$-tangle $T$ by a homeomorphism $\varphi : \partial(B,T_0)
\rightarrow \partial(B,T)$. We obtain a link $L=T_0 \cup T$ in 
a three-sphere $B {\cup}_{\varphi}B$ with 
a diagram $D(T_0 \cup T)$ as in Fig.~14. We may
assume that the diagram $D(T_0 \cup T)$ is connected. We color, 
in checkerboard fasion, the regions of the plane cut by the
diagram $D(T_0 \cup T)$ and  chose $n$ points $\{v_1,v_2,\cdots,v_n\}$ 
as in Fig.~15. Since $D(T_0 \cup
T)$ is connected, there is a spine $G$ of the 
black surface with the vertex set$\{v_1,v_2,\cdots,v_n\}$.
By retracting the black regions into the neighbourhood of $G$ 
and restricting to $B^3$, we have a 
required surface.

When we use the Seifert algorithm instead of checkerboard
coloring, we always obtain a bicollared disk-band representation. 
\hfill$\Box$

\begin{center}
\begin{tabular}{cc} 
\includegraphics[trim=0mm 0mm 0mm 0mm, width=.3\linewidth]
{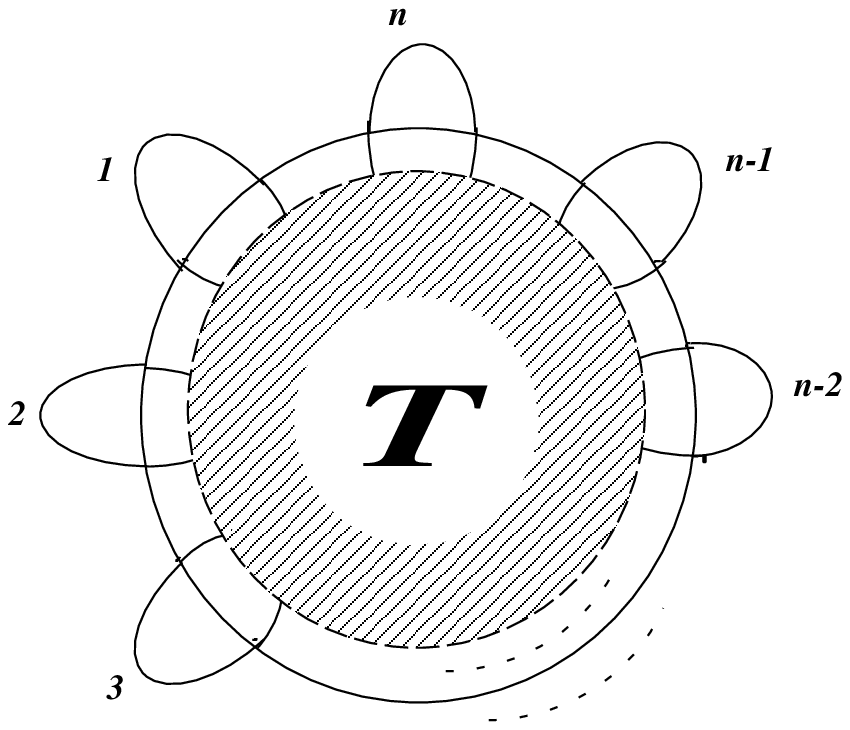}
&\hspace*{20mm}
\includegraphics[trim=0mm 0mm 0mm 0mm, width=.3\linewidth]
{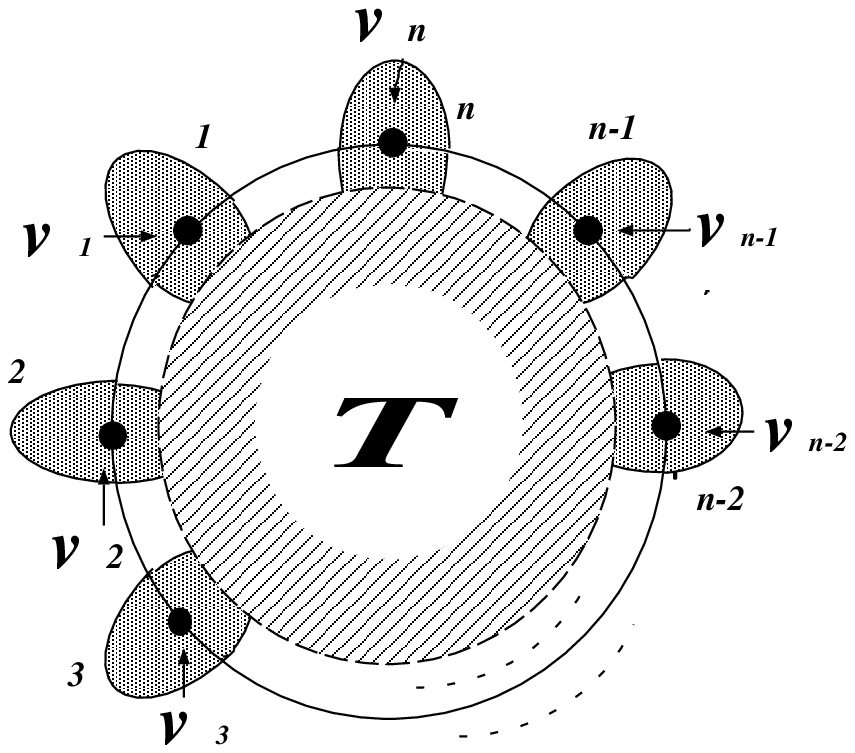}
\\
Fig. 14 &\hspace*{20mm} Fig. 15
\end{tabular}
\end{center}

\bigskip
\noindent
{\bf{2. Heegaard decompositions}}

\medskip
In addition to the surgery presentation, it is also useful to have another 
presentation of a $p$-fold cyclic branched cover. Our construction leads
straightforwardly to a Heegaard decomposition, 
a decomposition into a {\em compression body} \cite{B-O,C-G} 
and a handlebody, of a $p$-fold cyclic branched cover.

Let $F$ be a connected surface in $B^3$ bounded by an $n$-tangle 
$T$ and $n$ arcs in $\partial B^3$. The surface $F$ is defined to be 
{\em{free}} if the exterior of $F$ is homeomorphic to  
$(S_n\times I) \cup $(1-handles), where $S_n$ is an $n$-punctured sphere and 
 the all attaching points of the 1-handles are contained in 
 $S_n \times \{0\}$. As we observed before, any connected surface has 
 disk-band decomposition. A disk-band representation $\Omega (T_0;
\{D_1,\cdots,D_n\},\{b_1,\cdots,b_m\})$ is defined to be {\em{free}} if the surface
${\bigcup}_i D_i \cup {\bigcup}_j b_j $ is free.

First we consider the case of two-fold branched covers. 
The following algorithm gives a Heegaard
decomposition of $M^{(2)}(T)$. Start from a free disk-band
representation of $T$. Then we have a surgery description 
$\Sigma({\varphi}^{-1}({\bigcup}_i b_i), H_0)$ of
$M^{(2)}(T)$ (Fig.~3). 
Let $\alpha_1,\alpha_2,...,\alpha_{2m-n}$ be the $2m-n$ arcs  
which are the connected components of $T_0-\bigcup_j b_j$ 
contained in the interior of $B^3$. 
Then the complement of ${\varphi}^{-1}({\bigcup}_i b_i {\bigcup}_l
{\alpha}_l)$  and the arcs is a
compression body. Then $M^{(2)}(T) $ is obtained from the compression body 
by gluing the handle body as follows; (i) adding meridian disks of the arcs, 
and (ii) filling the rest according to the surgery description. 
For the tangle $T_1$ in Example 2, a Heegaard decomposition of
$M^{(2)}(T_1)$ is given as in Fig.~16. 

\begin{center}
\begin{tabular}{c} 
\includegraphics[trim=0mm 0mm 0mm 0mm, width=.45\linewidth]
{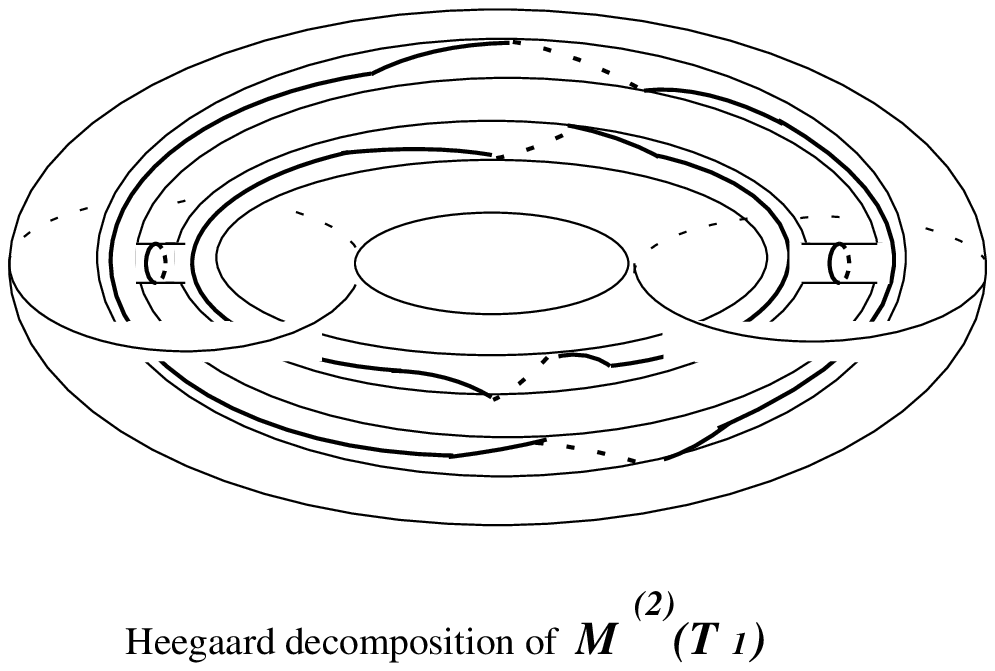}
\\
Fig. 16
\end{tabular}
\end{center}

\medskip
The similar method gives a Heegaard decomposition of a $p$-fold 
cyclic branched cover. 
Our construction is a modification of the construction in the
proof of Theorem 4. 
The handlebody part of the decomposition is
obtained from the handlebodies ${\bigcup}_j {\varphi}^{-1}(b_j \times
[-1,1]) $ by connecting them using $2m-n$ ``tubes'' 
along ${\varphi}^{-1}(T_0)$. We get genus $mp-n+1$ handlebody.

 From the observation above it follows that: 

\medskip
\noindent
{\bf{Theorem 7.}}
{\em Let $\Omega (T_0; \{D_1,\cdots,D_n\},\{b_1,\cdots,b_m\})$ be a 
free, bicollared disk-band representation of  an $n$-tangle $T$ in $B^3$.
Let $\varphi : H_0 \rightarrow B^3 $ be the $p$-fold cyclic 
branched cover of branched along $T_0$. 
Let $\alpha_1,\alpha_2,...,\alpha_{2m-n}$ be the $2m-n$ arcs  
which are the connected components of $T_0-\bigcup_j b_j$ 
contained in the interior of $B^3$. 
Then the following holds. \\
{\rm (a)} The complement $W=
\overline{H_0-\varphi^{-1}(\bigcup_j b_j\times [-1,1] \cup 
\bigcup_l N(\alpha_l))}$ is a compression body, where $N(\alpha_l)$ is 
the tubular neighbourhood of ${\alpha}_l$ in $B^3$.\\
{\rm (b)} The $p$-fold cyclic branched cover of $B^3$ branched 
along $T$ has a Heegaard decomposition into  
the compression body $W$ and a genus $mp-n+1$ handlebody. \\
{\rm (c)} The gluing map is given by the 
curves $c_{j,k}~(j=1,2,...,m,~k=1,2,...,p-1)$ and 
$m_l~(l=1,2,...,2m-n)$ in $\partial W$, where $c_{j,k}$ is the core of 
the annulus $b_{j,k}^+\cup b_{j,k+1}^-$ in Theorem 4 and 
$m_l$ is the meridian curve of $\varphi^{-1}(N(\alpha_l))$. 
\hfill$\Box$}  

\medskip
\noindent
In the theorem above, the assumption that a disk-band representation 
is bicollared is not necessary in the case that $p=2$.
The curves $c_{j,k}~(j=1,2,...,m,~k=1,2,...,p-1),~m_l~(l=1,2,...,2m-n)$ 
are essential, $mp-n+1$ of them are nonseparating and 
$m-1$ curves, $m_l$'s, are separating. 

\medskip
\noindent
{\bf{Remark 8.}} 
Since the surfaces given in the proof of Proposition 5 are 
connected and free, we can use them to find Heegaard decompositions of 
branched cyclic covers.
Let $c$ denote the crossing number of a connected diagram $D(T_0 \cup T), 
b$ the number of the black regions and $s$ the number of 
the Seifert circles of $D(T_0 \cup
T)$. Then we have a Heegaard decomposition of
$M^{(2)}(T)$ (resp. $M^{(p)}(T))$ of the genus $n+2c-2b+1$ 
(resp. $p(n+c-s)-n+1$).



\bigskip
{\small
\begin{thebibliography}{999}
\bibitem{A-K} S. Akbulut and R. Kirby: {\it{Branched covers of surfaces
	 in 4-manifolds}}, Math. Ann. 252 (1980), 111-131.

\bibitem{B-O} F. Bonahon and J. P. Otal: {\it{Scindements de Heegaard 
des espaces lenticulaires}}, Ann. Sci. Ec. Norm. Sup. 16 (1983), 451-466.

\bibitem{C-G} A. J. Casson and C. McA. Gordon: {\it{Reducing Heegaard
	splittings}}, Topology. Appl. 27 (1987), 275-283.

\bibitem{Con} J. H. Conway: {\it{An enumeration of knots and links and some 
	 of their related properties}}, 
Computational problems in Abstract Algebra, Proc. Conf. Oxford 1967,
Pergamon Press (1970), 329-358.

\bibitem {DJP} J. Dymara, T. Januszkiewicz, J. H. Przytycki: 
{\it{Symplectic structure on Colorings, Lagrangian 
tangles and Tits buildings }}, preprint (May 2001). 

\bibitem{Jac} W. Jaco: {\it{Lectures on three manifold topology}},
	 Conference board of Math. 43, A.M.S. (1980).

\bibitem{Mon} J. M. Montesinos: 
{\it{Surgery on links and double branched covers of $S^3$, 
Knots, groups and 3-manifolds}}, 
Ann. Math. Studies, Princeton Univ. Press. 84 (1975), 227-259.

\bibitem {P-S} V. V. Prasolov and A. B. Sossinsky: 
{\it{Knots, links, braids and 3-manifolds}}, A.M.S. (1997).

\bibitem {Rol} D. Rolfsen: 
{\it{Maps between 3-manifolds with nonzero degree: a new obstruction}}, 
in preparation for Proceedings of 
New Techniques in Topological Quantum Field Theory, 
NATO Advanced Research Workshop, 
August 2001, Canada. 
\end{thebibliography}}
\end{document}